\documentclass[12pt,a4paper]{article}
\usepackage{amsfonts}

\usepackage{amssymb}
\usepackage{color}

\usepackage{amsmath}

\newtheorem{theorem}{Theorem}[section]
\newtheorem{definition}[theorem]{Definition}
\newtheorem{lemma}[theorem]{Lemma}
\newtheorem{remark}[theorem]{Remark}
\newtheorem{corollary}[theorem]{Corollary}

\newcommand{\wt}{\operatorname{wt}}
\newcommand{\br}{\operatorname{br}}
\newcommand{\dd}{\operatorname{d}}
\newcommand{\lc}{\operatorname{lc}}
\usepackage{hyperref}

\begin{document}
\title{The Freiheitssatz and automorphisms for free brace algebras \footnote{
Supported by
the NNSF of China (Nos.  11501237; 11426112; 11401246),
the NSF of Guangdong Province of China (2016A030310099),
the Outstanding Young Teacher Training Program in Guangdong Universities (No. YQ2015155),
the Science and Technology Program of Huizhou City (2016X0429044,2017C0404020),
the Huizhou University (hzu201704, hzuxl201523).
}}
\author{Yu Li, Qiuhui Mo and Xiangui Zhao \\
{\small  Department of Mathematics, Huizhou
University}\\
{\small Huizhou 516007, P. R. China}\\
{\small Email: liyu820615@126.com }\\
{\small scnuhuashimomo@126.com}\\
{\small zhaoxg@hzu.edu.cn}}
 \date{}
\maketitle \noindent\textbf{Abstract:} 
Over a field of characteristic zero,
we prove that the Freiheitssatz holds for brace algebras,
the word problem for the brace algebras with a single defining relation is decidable,
two generated subalgebras of free brace algebras are free,
and that automorphisms of two generated free brace algebras are tame.

 \ \


\noindent {\bf AMS} Mathematics Subject Classification(2010): 
17A30, 
17A50, 
17A36, 
17D25 
\section{Introduction}
A brace algebra over a field is a vector space equipped with a family of linear operations satisfying some identities (see Definition \ref{def_Brace-Alg}).
Brace algebras have strong connections with other important classes of algebras.
For instance, brace algebras are used to prove Milnor-Moore type theorems for some Hopf algebras \cite{Ron01,Ron02};
a free brace algebra has also a free pre-Lie algebra structure (and thus has a free Lie structure) \cite{Foi10};
the pair of varieties (\emph{Brace}, \emph{Pre-Lie}) is a PBW-pair (in the sense of \cite{MS14}) \cite{LMZ18}.

In the present paper,
we continue the study of \cite{LMZ18} to investigate the Freiheitssatz (or independent theorem), word problem, subalgebras and automorphisms of brace algebras.

The Freiheitssatz for groups, one of the most important theorems of combinatorial group theory, was proved by Magnus \cite{M30} in 1930.
The Freiheitssatz states that: Let $G=gp\langle x_1,x_2,\dots,x_n\mid r=1\rangle$ be a group defined by a single cyclically reduced relator $r$.
If $x_n$ appears in $r$, then the subgroup of $G$ generated by $x_1,x_2,\dots,x_{n-1}$ is a free group with free generators  $x_1,x_2,\dots,x_{n-1}$.
As an application, the decidability of the word problem for one-relator groups was also proved by Magnus in the same paper.
After that, the Freiheitssatz for several other classes of algebras were established,
for example,
for Lie algebras (Shirshov \cite{Sh62b}),
for commutative (anti-commutative) algebras (Shirshov \cite{Sh62a}),
for associative algebras over a field of characteristic zero (proved by Makar-Limanov \cite{M-L85}, conjectured by Cohn \cite{Cohn62,Cohn74}),
for right-symmetric algebras (Kozybaev, Makar-Limanov and Umirbaev \cite{KMU08}),
for Poisson algebras and Novikov algebras over a field of characteristic zero (Makar-Limanov and Umirbaev \cite{MU11a,MU11b}),
and for generic Poisson algebras over a field of characteristic zero (Kolesnikov, Makar-Limanov and Shestakov \cite{KMS14}).
Using the Freiheitssatz and the decidability of the word problem for nonassociative algebras \cite{Z50}, Mikhalev and Shestakov \cite{ MS14} gave a uniform proof for the Freiheitssatz and the decidability of the word problem for commutative (anti-commutative) algebras, Akivis algebras and  Sabinin algebras.
Note that the Freiheitssatz for Poisson algebras in a positive characteristic is not true \cite{MU11a}.
The question about the decidability of the word problem for associative algebras (Poisson algebras, respectively) with a single defining relation and
the Freiheitssatz for associative algebras in a positive characteristic still remain open.

In Section \ref{sec_Freiheitssatz} of the present paper,
we prove the Freiheitssatz for brace algebras and the decidability of the word problem for the brace algebras with a single defining relation in characteristic zero.
These results {imply} the Freiheitssatz for right-symmetric algebras and the decidability of the word problem for the right-symmetric algebras {with a single defining relation} in characteristic zero \cite{KMU08}.

Recall that a variety of algebras is called {\it Schreier} if every subalgebra of a free algebra in this variety is also free.
It is known that the variety of pre-Lie algebras is not Schreier \cite{K07}.
Recently, Li, Mo and Zhao \cite{LMZ18} proved that the pair of varieties (\emph{Brace}, \emph{Pre-Lie}) over a field of characteristic zero
is a PBW-pair in the sense of \cite{MS14}.
Thus, by Theorem 1 of \cite{MS14},
the variety of brace algebras {in characteristic zero} is  not Schreier.
In a non-Schreier variety,
two generated subalgebras are not necessarily free.
For instance, two generated subalgebras of polynomial algebras
and associative algebras (both non-Schreier) are not necessarily free.
To our best knowledge,
the variety of right-symmetric algebras is the first non-Schreier one with the property that two generated subalgebras of a free algebra are free \cite{KMU08}.
We prove in Section \ref{sec_subalg} that
the variety of brace algebras also has this property,
or more precisely, two generated subalgebras of free brace algebras in characteristic zero are free.

It is well known \cite{Czer71-72,Jung42,Kulk53,M-L70} that the automorphisms of polynomial rings
and free associative algebras in two variables are tame.
Similar results concerning tame automorphisms are established for
free Poisson algebras and free right-symmetric algebras in two variables over a field of characteristic zero \cite{MTU09,KMU08}.
In the present paper,
we generalize this result to brace algebras:
automorphisms of two generated free brace algebras in characteristic zero are tame.

The rest of the paper is organized as follows.
In Section \ref{sec_brace-alg} we recall the definition and basic properties of brace algebras.
In Section \ref{sec_Freiheitssatz} we prove the Freiheitssatz and
study the word problem of brace algebras with a single defining relation.
In Section \ref{sec_subalg} we study subalgebras and automorphisms of brace algebras.
\section{Brace algebras and free brace algebras}
\label{sec_brace-alg}
\begin{definition}
\label{def_Brace-Alg}
\cite{Cha02,Ron01,Ron02}
A brace algebra is a couple $(A,\langle\rangle)$ where $A$ is a vector space and $\langle\rangle$ is a family of operators $A^{\otimes n}\longrightarrow
A$ defined for all $n\geq 2$:
$$
A^{\otimes n}\longrightarrow
A;a_n\otimes\dots\otimes a_1\longmapsto \langle a_n,\dots,a_{2};a_1\rangle,$$
with the following compatibilities: for all $a_1,\dots,a_m,b_1,\dots,b_n,c\in A$,
\begin{align*}
&\langle a_m,\dots,a_1;\langle b_n,\dots,b_1;c\rangle\rangle\\
=&\sum \langle V_{2n},\langle V_{2n-1};b_n\rangle,\dots,V_4,\langle V_3;b_2\rangle, V_2,\langle V_1;b_1\rangle,V_0;c\rangle,
\end{align*} where this sum runs over partitions of the ordered set $\{a_m,\dots,a_2,a_1\}$ into (possibly empty) consecutive intervals $V_{2n}\sqcup\dots\sqcup V_0$. We use the convention $\langle a\rangle=a$  for all $a\in A$.
\end{definition}

For example, if $A$ is a brace algebra and $a,b,c\in A$:
$$\langle a;\langle b;c\rangle \rangle=\langle a, b;c\rangle+\langle b,a;c\rangle+\langle \langle a;b\rangle;c \rangle.$$
From this relation, it follows immediately that
$(A, \langle -;-\rangle)$ is a left pre-Lie algebra.
Note that  a left (right) pre-Lie algebra is also called a left-symmetric (right-symmetric) algebra.
 Here is another example of relation in a brace algebra:
for all $a,b,c,d\in A$,
\begin{align*}
\langle a,b;\langle c;d\rangle\rangle=&\langle a,b,c;d\rangle+\langle a,\langle b;c\rangle;d\rangle+\langle \langle a,b;c\rangle;d\rangle\\
&+\langle a,c,b;d\rangle+\langle \langle a;c\rangle,b;d\rangle+\langle c,a,b;d\rangle.
\end{align*}

Let $X$ be a set.
Each letter $x\in X$ is called a {\it brace word of degree $1$}.
Let $u_i,1\leq i\leq n,$ be brace words of degrees $k_i,1\leq i\leq n$, respectively. Then
$w=\langle u_n,\dots,u_2;u_1\rangle$ is called {\it a brace word of degree $\dd(w):=\sum\limits_{i=1}^{n}k_i$}.
Denote by $\Omega(X)$ the set of all brace words on $X$.

Each letter $x$ in the alphabet $X$ is called a {\it normal brace word of degree $1$}.
Let $x\in X$, and $u_i,1\leq i\leq n,$ be normal brace words of degrees $k_i,1\leq i\leq n$ respectively. Then
$w=\langle u_n,\dots,u_2,u_1;x\rangle$ is called a {\it normal brace word of degree $1+\sum\limits_{i=1}^{n}k_i$}. Denote by $N(X)$ be the set of all normal brace words on $X$.

Let $F$ be a field of characteristic zero and $FN(X)$ the $F$-linear space spanned by $N(X)$. For normal brace words $v_i,1\leq i\leq m,1\leq m$ and $w=\langle u_n,\dots,u_2,u_1;x\rangle$, define
$$\langle v_m,\dots,v_2,v_1;w\rangle=\sum \langle V_{2n},\langle V_{2n-1};u_n\rangle,\dots,V_4,\langle V_3;u_2\rangle, V_2,\langle V_1;u_1\rangle,V_0;x\rangle,$$ where this
 sum runs over partitions of the ordered set $\{v_m,\dots,v_2,v_1\}$ into (possibly empty) consecutive intervals $V_{2n}\sqcup\dots\sqcup V_0$.
 Then $(FN(X),\langle\rangle)$ forms a free brace algebra generated by $X$ (see for example  \cite{Cha02,Foi02b,Foi10}), denoted by $Br(X)$.
 Each element of $Br(X)$ is also called a polynomial.

Let $X$ be a well ordered set and $w=\langle w_n,\dots,w_2,w_1;x\rangle$ a normal brace word on $X$.
Let $\br(w)=n+1$ be the \emph{breadth} of $w$. Then  we define
$$
\wt(w)=(\dd(w),\br(w),x,u_1,u_2,\dots,u_n)
$$
and order $N(X)$ by
$$
w>w' \Leftrightarrow \operatorname{wt}(w)>\wt(w') \ \ \mbox{lexicographically},
\text{ for all } w,w'\in N(X).
$$

This ordering is called degree
breadth inverse lexicographic ordering and used throughout this paper. It is easy to verify that this ordering is a well ordering.

For each nonzero polynomial $f\in Br(X)$, $f$ can be uniquely
presented as
$$
f=\alpha_1u_1+\alpha_2u_2+\dots+\alpha_nu_n,
$$
where $\alpha_i\in F$, $u_i\in N(X)$ for all $i$, $\alpha_1\neq
0$, $u_1>u_2>\dots>u_n$. Here, the
normal brace word $u_1$ is called the leading term of $f$, denoted
by $\bar{f}$ and $\alpha_1$ the leading coefficient of $f$, denoted
by $\lc(f)$.
The \emph{degree of $f$} is defined as the degree of its leading term,
i.e., $\dd(f):=\dd(u_1)$.
If $\lc(f)=1$, then  $f$ is called
a monic polynomial.

The rest of this section includes some elementary properties of the normal brace words and the degree
breadth inverse lexicographic ordering on $X$.

\begin{lemma}\label{l1'}
Let $w\in \Omega(X)$. Then, in
$Br(X)$, $w$ can be uniquely presented as
$$
w=a_1w_1+a_2w_2+\dots+a_n w_n,
$$
where $a_i$ is a positive integer and $w_i\in N(X)$ for each $i$.
\end{lemma}
\noindent\textbf{Proof.} Let us use induction on $\dd(w)$.
If $\dd(w)=1$, then $w\in X$ and the statement holds clearly. Let $w=\langle u_n,\dots,u_2,u_1;v\rangle$, where $v,u_1,\dots, u_n\in \Omega(X)$.
Obviously, the degree of each brace word belonging to $\{v,u_1,\dots,u_n\}$ is less than $\dd(w)$. Then by the inductive hypothesis we may assume without loss of generality that $v,u_1,\dots, u_n\in N(X)$.

If $v=x \in X$, then $w=\langle u_n,\dots,u_2,u_1;x\rangle\in N(X)$ and thus the statement holds.

If $v=\langle v_m,\dots,v_2,v_1;x\rangle$ with $x\in X$ and $m\geq 1$, then
\begin{align*}
w=&\langle u_n,\dots,u_2,u_1;\langle v_m,\dots,v_2,v_1;x\rangle\rangle\\
=&\sum \langle U_{2m},\langle U_{2m-1};v_m\rangle,\dots,U_4,\langle U_3;v_2\rangle, U_2,\langle U_1;v_1\rangle,U_0;x\rangle,
\end{align*}
where this sum runs over partitions of the ordered set $\{u_n,\dots,u_2,u_1\}$ into (possibly empty) consecutive intervals $U_{2n}\sqcup\dots\sqcup U_0$.
Then the statement follows from the inductive hypothesis immediately. \ \ \ $\square$

The following lemma appears in \cite{LMZ18} as Lemma 3.1.

\begin{lemma}\label{l1}
Let $v_i$ and $w=\langle u_n,\dots,u_2,u_1;x\rangle$ be normal brace words, $1\leq i\leq m$.
Then
$\overline{\langle v_m,\dots,v_2,v_1;w\rangle}$
is of the form $ \langle V_{n},u_n,\dots,V_2,u_2,\ V_1,u_1,V_0;x\rangle,$
where  $V_{n}\sqcup\dots\sqcup V_0$ is some (possibly empty) consecutive interval of the ordered set $\{v_m,\dots,v_2,v_1\}$.
\end{lemma}

\begin{remark}\label{r1} Lemma \ref{l1} can be also expressed as follows:
Let $v_i,1\leq i\leq m$ and $w=\langle u_n,\dots,u_2,u_1;x\rangle$ be normal brace words. Then
$\overline{\langle v_m,\dots,v_2,v_1;w\rangle}$
is of the form $ \langle U_{m},v_m,\dots,U_2,v_2,\ U_1,v_1, U_0;x\rangle,$
where  $U_{m}\sqcup\dots\sqcup U_0$ is some  (possibly empty) consecutive interval of the ordered set $\{u_n,\dots,u_2,u_1\}$.
\end{remark}

Immediately we have the following
\begin{corollary}\label{LMZ18-c1}\cite{LMZ18}
Let $v_i$ and $w$ be normal brace words for $1\leq i\leq m$. Then
$$
\br(\overline{\langle v_m,\dots,v_2,v_1;w\rangle})=m+\br(w).
$$
\end{corollary}

\begin{lemma}\label{l2}\cite{LMZ18}
Let $w,w',v_1,v_2,\dots,v_m,$ and $u$ be normal brace words. If $w>w'$, then
$$\overline{\langle v_m,\dots,v_1;w    \rangle}>\overline{\langle v_m,\dots,v_1;w'    \rangle} \text{ and}$$
$$\overline{\langle v_m,\dots,v_{i+1},w,v_i,\dots,v_1;u    \rangle}>\overline{\langle v_m,\dots,v_{i+1},w',v_i,\dots,v_1;u   \rangle}.$$
\end{lemma}

\begin{corollary}\label{LMZ18-c2}\cite{LMZ18}
Let $g, f_1,\dots,f_m\in Br(X)$. Then
$$
\overline{\langle f_m,\dots,f_2,f_1;g\rangle}=\overline{\langle \overline{f_m},\dots,\overline{f_2},\overline{f_1};\overline{g}\rangle}.
$$
\end{corollary}

\begin{lemma}\label{l3}
Let $u_i$, $v_j$ and $w$ be normal brace words for
$1\leq i\leq m$ and $1\leq j\leq n$.
 If $\overline{\langle u_m,\dots,u_1;w\rangle}=\overline{\langle v_n,\dots,v_1;w\rangle}$, then $m=n$ and $u_i=v_i$ for all $1\leq i\leq n$.
\end{lemma}
\noindent\textbf{Proof.}
Since $\overline{\langle u_m,\dots,u_1;w\rangle}=\overline{\langle v_n,\dots,v_1;w\rangle}$,
by Corollary \ref{LMZ18-c1} we have $m=n$ and $\sum_{i=1}^n{\dd(u_i)}=\sum_{i=1}^n{\dd(v_i)}$.
Let $w=\langle w_l,\dots,w_1;x\rangle$, each $w_i\in N(X)$, $0\leq l$ and $x\in X$.
 Then, by Remark \ref{r1}, $\overline{\langle v_n,\dots,v_1;w\rangle}$ is of the form $ \langle W_{n},v_n,\dots,W_2,v_2,\ W_1,v_1,W_0;x\rangle,$
where  $W_{n}\sqcup\dots\sqcup W_0$ is some consecutive interval of the ordered set $\{w_l,\dots,w_2,w_1\}$.
Clearly,  the polynomial $\langle u_n,\dots,u_1;w\rangle$ contains a  term $ \langle W_{n},u_n,\dots,W_2,u_2,\ W_1,u_1,W_0;x\rangle$ for the above $W_{0},\dots,W_n$.

Now, for a contradiction, suppose that there exists  some integer $t$ such that $u_1=v_1,\dots,u_{t-1}=v_{t-1}$ and $u_t> v_t$.
Then
$$\langle W_{n},u_n,\dots,W_2,u_2, W_1,u_1,W_0;x\rangle >\langle W_{n},v_n,\dots,W_2,v_2, W_1,v_1,W_0;x\rangle$$
since $\sum_{i=1}^n{\dd(u_i)}=\sum_{i=1}^n{\dd(v_i)}$, $u_1=v_1,\dots,u_{t-1}=v_{t-1}$ and $u_t> v_t$.
Thus $\overline{\langle u_m,\dots,u_1;w\rangle}> \langle W_{n},v_n,\dots,W_2,v_2, W_1,v_1,W_0;x\rangle=\overline{\langle v_n,\dots,v_1;w\rangle}$.
This is a contradiction.\ \ \ $\square$

\section{The Freiheitssatz for brace algebras}
\label{sec_Freiheitssatz}
In this section,
 over a field of characteristic zero,
we prove the decidability of the word problem for brace algebras with a single relation and the Freiheitssatz.

Let $X=\{x_1,\dots,x_M\}$ be a finite set and  we put $x_1<x_2<\dots<x_M$. Let $X_1=\{x_1,\dots,x_M,y\}$ and $x_M<y$.
For a brace word $w$ in the alphabet $X_1$, denote by $\dd_y(w)$ the degree of $w$ relative to $y$,
i.e., the number of occurrences of $y$ in $w$. Denote by $\Omega_y(X_1)$ the set of all brace words $u\in \Omega(X_1)$ with $d_y(u)=1$, and by $N_y(X_1)$ the set of all normal brace words $u\in N(X_1)$ with $\dd_y(u)=1$.

Let $f \in Br(X)$.  Define a brace algebra homomorphism $\psi: Br(X_1)\longrightarrow Br(X)$ by  $x_i\mapsto x_i (1\leq i\leq M), y\mapsto f$. Denote by $Id(f)$ the ideal generated by $f$ in $Br(X)$. Let $\psi (N_y(X_1))=B$.

\begin{lemma}\label{l4}
The ideal $Id(f)$ of $Br(X)$ is linearly spanned by $B$.
\end{lemma}
\noindent\textbf{Proof.} It is clear that $\psi(\Omega_y(X_1))$ linearly spans $Id(f)$.
By Lemma \ref{l1'}, each brace word $u \in \Omega_y(X_1)$ can be uniquely presented  as a linear combination of the elements in $N_y(X_1)$.
Thus the ideal
$Id(f)$ of $Br(X)$ is linearly spanned by $\psi (N_y(X_1))=B$.\ \ \ $\square$

\begin{lemma}\label{r2}
Let $u\in \Omega_y(X_1)$. Then $\psi(u)$ can be presented as a linear combination of $g_i\in B$ with $\overline{g_i}\leq \overline{\psi(u)}$.
\end{lemma}
\noindent\textbf{Proof.} It follows from  Lemma \ref{l1'} immediately.\ \ \ $\square$

Let $g,h\in Br(X)$ and $w\in N(X)$. Denote by $g\equiv h$ $\mod(B,w)$ if there exist  $\alpha,\beta, \gamma_i \in F$ and $g_i\in B$ such that $\alpha g-\beta h=\sum \gamma_i g_i$, where $\overline{g_i}<w$.

\begin{lemma}\label{l5}
If $g,h\in B$ and $\bar g=\bar h=w$, then $g\equiv h$ $\mod(B,w)$.
\end{lemma}
\noindent\textbf{Proof.} Let us use induction on $ w $. If $w=\bar f$, then $g=h=f$ and thus the statement holds obviously.

Suppose that $G,H\in N_y(X_1)$ and $\psi(G)=g,\psi(H)=h$.
Note that $G$ can be written in one of following forms:
$(i)$ $G=\langle G_m,\dots,G_1;y\rangle$; $(ii)$ $G=\langle G_m,\dots, G_{p+1},G_p,G_{p-1}\dots,G_1;x_i\rangle$, where only $G_p$ contains $y$.
Of course $H$ can be also written in one of following forms: $(i)$  $H=\langle H_n,\dots,H_1;y\rangle$; $(ii)$ $H=\langle H_n,\dots,H_{q+1},H_q,H_{q-1}\dots,H_1;x_j\rangle$,
where only $H_q$ contains $y$. Therefore we need to consider four cases.

Case 1. $G=\langle G_m,\dots,G_1;y\rangle$ and $H=\langle H_n,\dots,H_1;y\rangle$. In this case, $g=\psi(G)=\langle G_m,\dots,G_1;f\rangle$ and $h=\psi(H)=\langle H_n,\dots,H_1;f\rangle$. By Corollary \ref{LMZ18-c2}, $\bar g=\overline{\psi(G)}=\overline{\langle G_m,\dots,G_1;\bar f\rangle}$ and $\bar h=\overline{\psi(H)}=\overline{\langle H_n,\dots,H_1;\bar f\rangle}$.
Since $\bar g=\bar h$, by Lemma \ref{l3} we have $m=n$ and $G_s=H_s, 1\leq s\leq m$. Therefore $G=H$ and $g-h=0$.

Case 2. $G=\langle G_m,\dots, G_{p+1},G_p,G_{p-1}\dots,G_1;x_i\rangle$, where only $G_p$ contains $y$ and $H=\langle H_n,\dots,H_{q+1},H_q,H_{q-1},\dots,H_1;x_j\rangle$, where only $H_q$ contains $y$.
In this case, $g=\psi(G)=\langle G_m,\dots, G_{p+1},\psi(G_p),G_{p-1}\dots,G_1;x_i\rangle$ and $h=\psi(H)=\langle H_n,\dots,H_{q+1},\psi(H_q),H_{q-1},\dots,\\ H_1;x_j\rangle$.
By Corollary \ref{LMZ18-c2},
$\bar g=\langle G_m,\dots, G_{p+1},\overline{\psi(G_p)},G_{p-1}\dots,G_1;x_i\rangle$ and $\bar h=\langle H_n,\dots,\\ H_{q+1},\overline{\psi(H_q)}, H_{q-1}\dots,H_1;x_j\rangle$.
Since $\bar g=\bar h$, we have $x_i=x_j$ and $m=n$.

If $p=q$, then  $\bar g=\bar h$ implies that $\overline{\psi(G_p)}=\overline{\psi(H_p)}$ and $G_s=H_s, 1\leq s\leq m, s\neq p$. Since $\psi(G_p),\psi(H_p)\in B$ and  $\overline{\psi(G_p)}=\overline{\psi(H_p)}$, by the inductive hypothesis  we have $$\psi(G_p)\equiv \psi(H_p) \ \ \mod (B, \overline{\psi(G_p)}) .$$
If follows that
\begin{align*}
 & \langle G_m,\dots, G_{p+1},\psi(G_p),G_{p-1}\dots,G_1;x_i\rangle\\
\equiv& \langle G_m,\dots, G_{p+1},\psi(H_q),G_{p-1}\dots,G_1;x_i\rangle
\mod (B, \bar g).
\end{align*}
Since $G_s=H_s, 1\leq s\leq m, s\neq p$, we have
\begin{align*}
  &\langle G_m,\dots, G_{p+1},\psi(G_p),G_{p-1}\dots,G_1;x_i\rangle\\
\equiv& \langle H_m,\dots,H_{p+1},\psi(H_p),H_{p-1},\dots,H_1;x_i\rangle
\mod(B,\bar g),
\end{align*}
that is, $g\equiv h$ $\mod(B,\bar g)$.

Now consider the case $p\neq q$. Without loss of generality, we assume that $p>q$.
In this case, $\bar g=\bar h$ implies that $G_1=H_1,\dots,G_{q-1}=H_{q-1}$, $G_q=\overline{\psi(H_q)}$, $G_{q+1}=H_{q+1},\dots,G_{p-1}=H_{p-1}$, $\overline{\psi(G_p)}=H_p$, $G_{p+1}=H_{p+1},\dots,G_m=H_m$.
Since $G_q=\overline{\psi(H_q)}$ and $\overline{\psi(G_p)}=H_p$,
we may assume that $\psi(H_q)=\alpha G_q+\triangle$ and $\psi(G_p)=\beta H_p+\nabla$,
where $\alpha,\beta\in F$, $\triangle$ and $\nabla$ are linear combinations of normal brace words (on $X$) that are smaller than $G_q$ and $H_p$ respectively.
Therefore
\begin{align*}
\alpha g-\beta h=&\alpha g-\langle G_m,\dots, G_{p+1},\psi(G_p),G_{p-1},\dots,G_{q+1},\psi(H_q) ,G_{q-1},\dots,G_1;x_i\rangle\\
&+\langle G_m,\dots, G_{p+1},\psi(G_p),G_{p-1},\dots,G_{q+1},\psi(H_q) ,G_{q-1},\dots,G_1;x_i\rangle-\beta h\\
=&\langle G_m,\dots, G_{p+1},\psi(G_p),G_{p-1},\dots,G_{q+1},(\alpha G_q-\psi(H_q)) ,G_{q-1},\dots,G_1;x_i\rangle\\
&+\langle G_m,\dots, G_{p+1},(\psi(G_p)-\beta H_p),G_{p-1},\dots,G_{q+1},\psi(H_q) ,G_{q-1},\dots,G_1;x_i\rangle\\
=&\langle G_m,\dots, G_{p+1},\psi(G_p),G_{p-1},\dots,G_{q+1},\triangle ,G_{q-1},\dots,G_1;x_i\rangle\\
&+\langle H_m,\dots, H_{p+1},\nabla,H_{p-1},\dots,H_{q+1},\psi(H_q) ,H_{q-1},\dots,H_1;x_i\rangle.
\end{align*}
Hence $\alpha g- \beta h$  is a linear combination of elements $g_l\in B$ where all $\bar {g_l}<\bar g=w$, since $\overline{\triangle}<G_q$ and $\overline{\nabla}<H_p$.
Therefore $g\equiv h$ $\mod(B,w)$.

Case 3. $G=\langle G_m,\dots,G_1;y\rangle$, and $H=\langle H_n,\dots,H_{q+1},H_q,H_{q-1},\dots,H_1;x_j\rangle$, where only $H_q$ contains $y$. In this case $g=\psi(G)=\langle G_m,\dots,G_1;f\rangle$ and
$h=\psi(H)=\langle H_n,\dots,H_{q+1},\psi(H_q),H_{q-1},\dots,H_1;x_j\rangle$. Suppose that $f=\alpha_{\bar f}\bar f  +\delta$, where $\bar f=\langle u_t,\dots, u_2,\\
u_1;x_i\rangle$ and $\bar \delta<\bar f$.
By Remark \ref{r1}, $\bar g$ is of the form $\langle U_m,G_m,\dots,U_1,G_1,U_0;x_i\rangle$,
where  $U_{m}\sqcup\dots\sqcup U_0$ is some consecutive interval of the ordered set $\{u_t,\dots,u_2,u_1\}$.
Then we have that $x_i=x_j$ and $G_s=\overline{\psi(H_q)}$ for some $s$, since $\bar g=\bar h$ and $\bar f \leq \overline{\psi(H_q)}$ . Let us assume that $\psi(H_q)=\alpha G_s+\triangle$, where $\triangle$ is a linear combination of normal brace words on $X$ which are smaller than $G_s$.
Then
\begin{align*}
\alpha g=&\langle G_m,\dots, G_{s+1},\alpha G_{s},G_{s-1},\dots, G_1;f\rangle\\
=&\langle G_m,\dots, G_{s+1},(\alpha G_{s}+\triangle),G_{s-1},\dots, G_1;f\rangle-\langle G_m,\dots, G_{s+1},\triangle,G_{s-1},\dots, G_1;f\rangle\\
=&\langle G_m,\dots, G_{s+1},\psi(H_q),G_{s-1},\dots, G_1;\alpha_{\bar f} \bar f\rangle+\langle G_m,\dots, G_{s+1},\psi(H_q),G_{s-1},\dots, G_1;\delta\rangle\\
&-\langle G_m,\dots, G_{s+1},\triangle,G_{s-1},\dots, G_1;f\rangle.\\
\end{align*}
Since $\bar \delta<\bar f$, by Lemma \ref{r2},
 $\langle G_m,\dots, G_{s+1},\psi(H_q),G_{s-1},\dots, G_1;\delta\rangle$ can be presented as a linear combination of $g_t\in B$, where $\overline {g_t}\leq \overline{\langle G_m,\dots, G_{s+1},\psi(H_q),G_{s-1},\dots, G_1;\delta\rangle}<\bar g$.
Clearly, $\langle G_m,\dots, G_{s+1},\triangle,G_{s-1},\dots, G_1;f\rangle$ can be also presented as a linear combination of $g_t'\in B$ with $\overline {g_t'}<\bar g$ since $\overline{\triangle}<G_s$.
By Lemma \ref{r2},
\begin{align*}
&\langle G_m,\dots, G_{s+1},\psi(H_q),G_{s-1},\dots, G_1;\alpha_{\bar f} \bar f\rangle\\
=&\alpha_{\bar f} \sum \langle V_m,G_m,\dots, V_{s+1},G_{s+1},V_{s},\psi(H_q),V_{s-1},G_{s-1},\dots, V_1,G_1,V_0;x_i\rangle+\sum_l\alpha_lg_l,
\end{align*}
where the first sum runs over partitions of the ordered set $\{u_t,\dots,u_2,u_1\}$ into (possibly empty) consecutive intervals $V_{m}\sqcup\dots\sqcup V_0$ and $g_i\in B$, $\overline{g_l}<\bar g $.
It is easy to see that for each consecutive intervals $V_{m}\sqcup\dots\sqcup V_0$
\begin{align*}
&\overline{\langle V_m,G_m,\dots, V_{s+1},G_{s+1},V_{s},\psi(H_q),V_{s-1},G_{s-1},\dots, V_1,G_1,V_0;x_i\rangle}\\
\leq& \overline{\langle G_m,\dots, G_{s+1},\psi(H_q),G_{s-1},\dots, G_1;\alpha_{\bar f} \bar f\rangle}=\bar g=\bar h.
\end{align*}
Therefore by the same argument as in Case 2 , we have
$
\alpha g=\beta h+\gamma_rg_r,
$
where $\beta,\gamma_r\in F$,  $g_r\in B$ and $\overline{g_r}<\bar h=\bar g=w$.
Therefore $g\equiv h$ $\mod(B,w)$.

Case 4. $G=\langle G_m,\dots, G_{p+1},G_p,G_{p-1}\dots,G_1;x_i\rangle$, where only $G_p$ contains $y$, and $H=\langle H_n,\dots,H_1;y\rangle$. The statement in this case can be proved similarly as in Case 3.\ \ \ $\square$

\begin{lemma}\label{l5'}
Let $h\in Id(f)$ and $h\neq 0$. Then there exists some $g\in B$ such that $\bar h=\bar g$.
\end{lemma}
\noindent\textbf{Proof.} Let $h\in Id(f)$ and $h\neq 0$. We may assume, by Lemma \ref{l4}, that $h=\sum_{i=1}^n\alpha_ig_i$, where $\alpha_i\in F$ and $g_i\in B$.
Suppose that $\overline{g_1}=\overline{g_2}=\dots=\overline{g_l}> \overline{g_{l+1}} \ge\dots$. Let us use induction on $\overline {g_1}$.

If $l=1$, then $\bar h=\overline{g_1}$ and hence the statement holds.

If $l>1$, then by Lemma \ref{l5}, $g_j=\beta_j g_1+\sum_{k_j}\gamma_{k
_ j}h_{k_ j}$,  where $j=2,\dots,l$, $\beta_j, \gamma_{k_j}\in F$, $h_{k_ j}\in B$ and $\overline{h_{k_ j}}<\overline {g_1}$. Therefore $h=(\alpha_1+\sum_{i=2}^l\alpha_i\beta_i)g_1+\triangle$, where $\triangle$ is a linear combination of $g_j'\in B$ and $\overline{g_j'}<\overline{g_1}$. If $\alpha_1+\sum_{i=2}^l\alpha_i\beta_i\neq 0$, then $\bar h=\overline{g_1}$. If $\alpha_1+\sum_{i=2}^l\alpha_i\beta_i= 0$, then the statement follows from the inductive hypothesis.\ \ \ $\square$

From the above lemma, we immediately have the following
\begin{corollary}\label{c1}
Let $h\in Br(X)$ and $h\neq 0$. If $\dd(h)<\dd(f)$, then $h\notin Id(f)$.
\end{corollary}

L.A. Bokut \cite{Bo72} proved the
undecidability of the word problem for Lie algebras.
An explicit example of a finitely presented Lie algebra with the
undecidable word problem was constructed by G. P. Kukin \cite{Kukin77} (see also \cite{CLT17}).
The undecidability word problem for right-symmetric algebras follows directly from Segal's  analogue of the
Poincar$\acute{e}$-Birkhoff-Witt theorem for right-symmetric algebras  \cite{Se94} and Bokut's result.
Li, Mo and Zhao \cite{LMZ18} proved that the pair of varieties (\emph{Brace}, \emph{Pre-Lie}) is a PBW-pair in the sense of \cite{ MS14}.
Together with the undecidability word for right-symmetric algebras, it follows that the word problem for brace algebras is also undecidable.
On the other hand, Shirshov \cite {Sh62a} proved the decidability of the word problem for Lie algebras with a single defining relation. Kozybaev, Makar-Limanov and Umirbaev \cite{KMU08} proved the decidability of the word problem for right-symmetric algebras with a single defining relation.
In the case of brace algebras over a field of characteristic zero, we have the following result.

\begin{theorem}
The word problem for brace algebras with a single defining relation is decidable.
\end{theorem}
\noindent\textbf{Proof.} Let $h\in Br(X)$ and $h\neq 0$. If  $\dd(h)< \dd(f)$, then, by Corollary \ref{c1}, $h\notin Id(f)$. Now we assume that $\dd(h)\ge \dd(f)$. Obviously, there exist only finitely many elements $g\in B$ such that $\dd(h)=\dd(g)$. Hence we can effectively determine whether there exists some element $g\in B$ such that $\bar h= \bar g$.
If there does not exist $g\in B$ such that $\bar h= \bar g$, then by Lemma \ref{l5'}, $h\notin Id(f)$. If there exists some $g\in B$ such that $\bar h= \bar g$, then let $h_1=h-\alpha g$ where $\alpha \in F$ satisfying
 $\lc(h)=\alpha \lc(g)$. Clearly, we have that $\overline{h_1}<\bar h$, and $h\in Id(f)$ if and only if $h_1\in Id(f)$.
 Note that $<$ is a well ordering.
 Therefore we can effectively determine whether $h\in Id(f)$.\ \ \ $\square$

\begin{lemma}\label{l6}
Given $u\in N(X)$, let $\phi_u: Br(X)\rightarrow Br(X)$ be a brace homomorphism defined by $x_i\mapsto x_i (1\leq i<M), x_M\mapsto \langle u;x_M\rangle$. If $v,w\in N(X)$ and $\overline{\phi_u(v)}= \overline{\phi_u(w)}$, then $v=w$.
\end{lemma}
\noindent\textbf{Proof.} If $v=x_i$, then $\overline{\phi_u(v)}= \overline{\phi_u(w)}$ implies that $w=x_i$, and hence $v=w$.
Let $v=\langle v_s,\dots, v_1; x_i\rangle$ and $w=\langle w_t,\dots, w_1; x_j\rangle$. By Corollary \ref{LMZ18-c2}, $\overline{\phi_u (v)}=\overline{\langle \overline{\phi_u (v_s)},\dots, \overline{\phi_u(v_1)}; \overline{\phi_u(x_i)}\rangle}$
and $\overline{\phi_u (w)}=\overline{\langle \overline{\phi_u (w_t)},\dots, \overline{\phi_u(w_1)}; \overline{\phi_u(x_j)}\rangle}$.

If $1\le i<M$, then $\overline{\phi_u (v)}=\langle \overline{\phi_u (v_s)},\dots, \overline{\phi_u(v_1)} ;x_i\rangle$. From $\overline{\phi_u(v)}= \overline{\phi_u(w)}$ and Lemma \ref{l3}, it follows that $i =j$, $s=t$ and
$\overline{\phi_u(v_1)}=\overline{\phi_u(w_1)},\dots,\overline{\phi_u (v_s)}=\overline{\phi_u (w_s)}$. Then by the inductive hypothesis  on $\dd(v)$,  we have $v_1=w_1,\dots,v_s=w_s$, and hence $v=w$.

If $i=M$, then obviously we have $j=M, s=t$ since $\overline{\phi_u(v)}= \overline{\phi_u(w)}$. Then by Lemma \ref{l3} $\overline{\phi_u(v_1)}=\overline{\phi_u(w_1)},\dots,\overline{\phi_u (v_s)}=\overline{\phi_u (w_s)}$. Then
by the inductive hypothesis  on $\dd(v)$ again we get $v_1=w_1,\dots,v_s=w_s$, and hence $v=w$.\ \ \ $\square$

As we mentioned in the introduction,
Shirshov \cite {Sh62a} proved the Freiheitssatz for Lie algebras. Kozybaev, Makar-Limanov and Umirbaev \cite{KMU08} proved the Freiheitssatz for right-symmetric algebras.
In the case of brace algebras over a field of characteristic zero, we have the following result.

\begin{theorem} (Freiheitssatz)
Let $Br(x_1,x_2,\dots, x_M)$ be the free brace algebra over a field $F$ of characteristic 0 in the variables $\{x_1,x_2,\dots, x_M\}$.
If $f\in Br(x_1,x_2,\dots,x_{M})$ and $f\notin Br(x_1,x_2,\dots,x_{M-1})$, then $Id(f)\cap Br(x_1,x_2,\dots,x_{M-1})=0$.
\end{theorem}
\noindent\textbf{Proof.}
Suppose that there exists some $h\in Br(X), h\neq 0$ such that  $h\in Id(f)\cap Br(x_1,x_2,\dots,x_{M-1})$.
We choose a normal brace word $u\in N(X)$ with $\dd(u)\ge \dd(h)$, and  then consider the endomorphism $\phi_u$ of $Br(X)$  defined by $x_i\mapsto x_i (1\leq i<M), x_M\mapsto \langle u;x_M\rangle$. It is clear that $\phi_u(h)=h\in Id(\phi_u(f))\cap Br(x_1,x_2,\dots,x_{M-1})$. By Lemma \ref{l6}, we have $\dd(\phi_u(f))>\dd(h)$. Therefore by Corollary \ref{c1},  $h\notin Id(\phi_u(f))$. This is a contradiction.\ \ \ $\square$

The next is a direct formulation of the Freiheitssatz for brace algebras in the language of freeness.
\begin{corollary} (Freiheitssatz)
Let $Br(x_1,x_2,\dots, x_M)$ be the free brace algebra over a field $F$ of characteristic 0 in the variables $\{x_1,x_2,\dots, x_M\}$.
If $f\in Br(x_1,x_2,\dots,x_{M})$ and $f\notin Br(x_1,x_2,\dots,x_{M-1})$, then the subalgebra of the quotient algebra
$Br(x_1,x_2,\dots, x_M)/Id(f)$ generated by $x_1+Id(f),x_2+Id(f),\dots,x_{M-1}+Id(f)$ is a free brace algebra with free generators $x_1+Id(f),x_2+Id(f),\dots,x_{M-1}+Id(f)$.
\end{corollary}

\section{Subalgebras and automorphisms of free brace algebras}
\label{sec_subalg}
Remember that a variety of algebras is called {\it Schreier} if every subalgebra
of a free algebra in this variety is also free.
Kozybaev proved in \cite{K07} that the variety of pre-Lie algebras is not a Schreier variety.
Li, Mo and Zhao \cite{LMZ18} proved that the pair of varieties (\emph{Brace}, \emph{Pre-Lie}) is a PBW-pair in the sense of \cite{ MS14}.
Then, by the Theorem 1 of \cite{MS14}, we know that the variety of brace algebras in characteristic zero is  not a Schreier variety.
However, we prove in this section that two generated subalgebras of free brace algebras in characteristic zero  are free. We also prove that automorphisms of two generated free brace algebras in characteristic zero are tame.

A subset $S$ of a brace algebra is called {\it algebraically independent}
if the elements of $S$ do not satisfy any non-trivial brace polynomial equation.

Let $A$ be a brace algebra and $S$ a subset of $A$. Denote by $alg_{A}(S)$ the subalgebra of $A$ generated by $S$.
Let $X=\{x_1,x_2,\dots, x_M\}$.
\begin{lemma}\label{l7}
Let $f\in Br(X)$ be a non-zero element. Then $alg_{Br(X)}(f)$ is a free brace algebra with a free generator $f$.
\end{lemma}
\noindent\textbf{Proof.} Assume that $\{f\}$ is algebraically dependent. Then there exists a non-zero element $p(y)=\alpha_1W_1(y)+\alpha_2W_2(y)+\dots+\alpha_nW_n(y) $ of  $Br(y)$, where each $W_l(y)$ is a normal brace word on $\{y\}$, such that $p(f)=0.$ It follows that there are two different normal brace words $W_i(y)$ and $W_j(y)$ $(i\neq j)$ such that $\overline{W_i(f)}=\overline{W_j(f)}$.
Let us assume that $w_a(y)$ and $w_b(y)$ is a pair of different normal brace words on $y$ with this property and the minimal degree $\dd(w_a(y))+\dd(w_b(y))$.
We can write $w_a(y)=\langle w_s(y),\dots,w_1(y);y\rangle$ and $w_b(y)=\langle w'_t(y),\dots,w'_1(y);y\rangle$, where $w_i(y)$ and $w'_j(y)$ are normal brace words on $\{y\}$ for all $i,j,(1\leq i\leq s, 1\leq j\leq t)$. By Corollary \ref{LMZ18-c2}, $\overline{w_a(f)}=\overline{w_b(f)}$
implies that $\overline{\langle \overline{w_s(f)},\dots,\overline{w_1(f)};\overline{f}\rangle}=\overline{\langle \overline{w'_t(f)},\dots,\overline{w'_1(f)};\overline{f}\rangle}$. Then according to Lemma \ref{l3} we have $s=t$ and $\overline{w_i(f)}=\overline{w'_i(f)}$ for each $i,(1\leq i\leq s)$. Since $\dd(w_i(y))+\dd(w'_i(y))<\dd(w_a(y))+\dd(w_b(y))$, we may conclude that  $w_i(y)=w'_i(y)$ for each $i,(1\leq i\leq s)$, and then  $w_a(y)=w_b(y)$. This is a contradiction.\ \ \ $\square$

Given two non-zero elements $f_1,f_2\in Br(X)$, let $\psi$ be a brace algebra homomorphism from $Br(y_1,y_2)$ to $Br(X)$ defined by  $y_i\mapsto f_i$ $(1\leq i\leq 2)$.

\begin{lemma}\label{l8}
If there exist two different normal brace words $W_i$ and $W_j$ on $\{y_1,y_2\}$  such that $\overline{\psi (W_i)}=\overline{\psi (W_j)}$,
then there exists a normal brace word $q(y)$ on $\{y\}$ such that $\overline{f_1}=\overline{q(f_2)}$ or $\overline{f_2}=\overline{q(f_1)}$.
\end{lemma}
\noindent\textbf{Proof.} If $\overline{f_1}=\overline{f_2}$, then the statement holds clearly.

Now, we assume that $\overline{f_1}>\overline{f_2}$. Let $w_a$ and $w_b$ be a pair of different normal brace words on $\{y_1,y_2\}$ with the property $\overline{\psi(w_a)}=\overline{\psi(w_b)}$ and the minimal degree $\dd(w_a)+\dd(w_b)$.

Without loss of generality, we may assume that $w_a=\langle w_s,\dots,w_1;y_1\rangle$, where each $w_i,(1\leq i\leq s)$ is a normal brace word on
$\{y_1,y_2\}$. Note that $w_b$ has two possible forms: (i) $w_b=\langle w'_t,\dots,w'_1;y_1\rangle$; (ii) $w_b=\langle w'_t,\dots,w'_1;y_2\rangle$.

If $w_b=\langle w'_t,\dots,w'_1;y_1\rangle$, then by $\overline{\psi (w_a)}=\overline{\psi (w_b)}$ and Corollary \ref{LMZ18-c2}
$$\overline{\langle \overline{\psi (w_s)},\dots,\overline{\psi (w_1)};\overline{f_1}\rangle}=\overline{\langle \overline{\psi (w'_t)},\dots,\overline{\psi (w'_1)};\overline{f_1}\rangle}.$$
According to Lemma \ref{l3} we have $s=t$ and $\overline{\psi (w_i)}=\overline{\psi (w'_i)}$ for each $i,(1\leq i\leq s)$.
Since $\dd(w_i)+\dd(w'_i)< \dd(w_a)+\dd(w_b)$ for each $i,(1 \leq i\leq s)$ and $w_a, w_b$ is a minimal pair, we may conclude that $w_i=w'_i$ for each $i,(1\leq i\leq s)$ and thus $w_a=w_b$. This is a contradiction.

Let us assume that $w_b=\langle w'_t,\dots,w'_1;y_2\rangle$. Then by $\overline{\psi (w_a)}=\overline{\psi (w_b)}$ and Corollary \ref{LMZ18-c2}
$$\overline{\langle \overline{\psi (w_s)},\dots,\overline{\psi (w_1)};\overline{f_1}\rangle}=\overline{\langle \overline{\psi (w'_t)},\dots,\overline{\psi (w'_1)};\overline{f_2}\rangle}.$$
Suppose that $\overline{f_1}=\langle u_p,\dots,u_1;x_i\rangle$ and $\overline{f_2}=\langle v_q,\dots,v_1;x_j\rangle$,
where $x_i,x_j\in X$,  $u_s, (1\leq s\leq p)$ and $v_t, (1 \leq t\leq q)$ are normal brace words on $X$.

Since $\overline{f_1}>\overline{f_2}$, we should have $\overline{\psi (w_s)},\dots,\overline{\psi (w_1)},\overline{\psi (w'_t)},\dots,\overline{\psi (w'_1)}\geq \overline{f_2}>v_k, (1\leq k\leq q)$.
It follows that $\overline{\langle \overline{\psi (w'_t)},\dots,\overline{\psi (w'_1)};\overline{f_2}\rangle}=
\langle v_q,\dots,v_1,\overline{\psi (w'_t)},\dots,\overline{\psi (w'_1)};x_j\rangle$.
By Remark \ref{r1}, we may assume that
$$
\overline{\langle \overline{\psi (w_s)},\dots,\overline{\psi (w_1)};\overline{f_1}\rangle}= \langle U_{s},\overline{\psi (w_s)},\dots,U_2,\overline{\psi (w_2)},\ U_1,\overline{\psi (w_1)}, U_0;x_i\rangle,
$$
where  $U_{s}\sqcup\dots\sqcup U_0$ is some  (possibly empty) consecutive interval of the ordered set $\{u_p,\dots,u_2,u_1\}$.

Since $\overline{\langle \overline{\psi (w_s)},\dots,\overline{\psi (w_1)};\overline{f_1}\rangle}=\overline{\langle \overline{\psi (w'_t)},\dots,\overline{\psi (w'_1)};\overline{f_2}\rangle}$ and  $\overline{\psi (w_s)},\dots,\overline{\psi (w_1)}>v_k, (1\leq k\leq q)$,
we may conclude that $x_i=x_j$, $s\leq t$ and $\overline{\psi (w_1)}=\overline{\psi (w'_{l_1})},\dots, \overline{\psi (w_s)}=\overline{\psi (w'_{l_s})}$,
for some integers ${l_1},\dots,{l_s}$, where $1\leq {l_1}<\dots<{l_s}\leq t$. Then we have
\begin{align*}
\overline{f_1}=&\langle U_{s},\dots,U_2,\ U_1,U_0;x_i\rangle\\
=&\langle v_q,\dots,v_1,\overline{\psi (w'_t)},\dots,\overline{\psi (w'_{l_s+1})},\overline{\psi (w'_{l_s-1})},\dots,\overline{\psi (w'_{l_1+1})},\overline{\psi (w'_{l_1-1})},\dots,\overline{\psi (w'_1)};x_i\rangle\\
=&\overline{\langle\overline{\psi (w'_t)},\dots,\overline{\psi (w'_{l_s+1})},\overline{\psi (w'_{l_s-1})},\dots,\overline{\psi (w'_{l_1+1})},\overline{\psi (w'_{l_1-1})},\dots,\overline{\psi (w'_1)};\langle v_q,\dots,v_1;x_i\rangle\rangle}\\
=&\overline{\langle\overline{\psi (w'_t)},\dots,\overline{\psi (w'_{l_s+1})},\overline{\psi (w'_{l_s-1})},\dots,\overline{\psi (w'_{l_1+1})},\overline{\psi (w'_{l_1-1})},\dots,\overline{\psi (w'_1)};\overline{f_2}\rangle}.
\end{align*}
From the above equalities, it follows that $w'_t,\dots,w'_{l_s+1},w'_{l_s-1},\dots,w'_{l_1+1},w'_{l_1-1},\dots,w'_1$ are normal brace words on $\{y_2\}$.
Let $q(y_2)=\langle w'_t,\dots,w'_{l_s+1},w'_{l_s-1},\dots,w'_{l_1+1},w'_{l_1-1},\dots,w'_1;y_2\rangle$. clearly, $\overline{f_1}=\overline{q(f_2)}$.

The statement in the case $\overline{f_1}<\overline{f_2}$ can be proved in a similar way.\ \ \ $\square$

\begin{theorem}\label{t4.2}
Let $f_1,f_2\in Br(X)$ be two non-zero elements. Then $alg_{Br(X)}(f_1,f_2)$ is a free brace algebra.
\end{theorem}
\noindent\textbf{Proof.} If $\{f_1,f_2\}$ is algebraically independent, then the statement holds clearly.

Let us assume that $\{f_1,f_2\}$ is algebraically dependent.
Then there exists a non-zero element $p(y_1,y_2)=\alpha_1W_1+\alpha_2W_2+\dots+\alpha_mW_m$,
where each $W_l$ is a normal brace word on $\{y_1,y_2\}$, such that $p(f_1,f_2)=0$. This implies that there are two different normal brace words $W_i$ and $W_j$ such that $\overline{W_i(f_1,f_2)}=\overline{W_j(f_1,f_2)}$. By Lemma \ref{l8}, there exists a normal brace word $q(y)$ on $\{y\}$ such that $\overline{f_1}=\overline{q(f_2)}$ or $\overline{f_2}=\overline{q(f_1)}$.

If $\overline{f_1}=\overline{q(f_2)}$, then set $g:=f_1-\alpha q(f_2)$, where $\alpha\in F$ and $lc(f_1)=\alpha lc(q(f_2))$.
Clearly, $g<\overline{f_1}$, $alg_{Br(X)}(f_1,f_2)=alg_{Br(X)}(g,f_2)$ and $\{g,f_2\}$ is also algebraically dependent.

For the other case, we set $g:=f_2-\alpha q(f_1)$, where $\alpha\in F$ and $lc(f_2)=\alpha lc(q(f_1))$.

Then after finite times of substitution on the generators of $alg_{Br(X)}(f_1,f_2)$, we have $alg_{Br(X)}(f_1,f_2)=alg_{Br(X)}(0,h)$. Therefore, by Lemma \ref{l7}, $alg_{Br(X)}(f_1,f_2)$ is a free brace algebra with a free generator $h$. \ \ \ $\square$

Let $\mathcal{M}=\{f_1,\dots,f_m\}$ be a subset of $Br(X)$. Then the transformation  $$f_j\mapsto f_j,j\neq {i}, f_{i}\mapsto \alpha f_{i}+g(f_{1},\dots,f_{i-1},f_{i+1},\dots,f_m)$$
where $0\neq\alpha\in F$ and $g\in Br(y_1,\dots,y_{i-1},y_{i+1},\dots,y_m)$, is called an elementary transformation of $\mathcal{M}$.

Recall that an automorphism $\phi$ of a brace algebra $Br(x_1,\dots,x_M)$ is called {\it elementary} if $\phi(x_j)=x_j$ for any $j\neq i$ and $\phi(x_i)=\alpha x_i+f$, where $f\in Br(x_1,\dots,x_{i-1},x_{i+1},\dots,x_M)$. Automorphisms which can be expressed as a composition of elementary automorphisms are called {\it tame}. Non-tame automorphisms are called {\it wild}.

Denote by $\phi=(f_1,f_2,\dots,f_M)$ the automorphism of $Br(x_1,\dots,x_M)$ defined by $\phi(x_i)=f_i, 1\leq i\leq M$.
It is well known that $\phi$ is tame if and only if there exists a finite sequence of elementary transformations such that
$$
(f_1,f_2,\dots,f_M)\rightarrow \dots \rightarrow (x_1,\dots,x_M).
$$

\begin{theorem}
Automorphisms of two generated free brace algebras are tame.
\end{theorem}
\noindent\textbf{Proof.} Let $\varphi=(f_1,f_2)$ be an automorphism of $Br(x_1,x_2)$. If there exists some normal brace word $q(y)
$ on $\{y\}$ such that $\overline{f_1}=\overline{q(f_2)}$ or $\overline{f_2}=\overline{q(f_1)}$, then we use elementary transformation $(f_1,f_2)\rightarrow (f_1-lc(f_1)lc(q(f_2))^{-1}q(f_2),f_2)$ or
$(f_1,f_2)\rightarrow (f_1,f_2-lc(f_2)lc(q(f_1))^{-1}q(f_1))$. After a finite number of elementary transformations we obtain a pair $g_1,g_2$, where $\overline{g_1}\neq \overline{q(g_2)}$ and $\overline{g_2}\neq \overline{q(g_1)}$ for any normal brace word $q(y)$ on $\{y\}$. Without loss of generality, we assume here that $\overline{g_1}<\overline{g_2}$.
Clearly, $Br(x_1,x_2)$ is also generated by $g_1,g_2$. Therefore $x_1=p_1(g_1,g_2)$ and $x_2=p_2(g_1,g_2)$ for brace polynomials $p_1(y_1,y_2), p_2(y_1,y_2)\in Br(y_1,y_2)$, and thus $x_1=\overline{x_1}=\overline{p_1(g_1,g_2)}$ and $x_2=\overline{x_2}=\overline{p_2(g_1,g_2)}$. By Lemma \ref{l8}, we know that distinct normal brace words on $\{g_1,g_2\}$ have distinct leading terms in $Br(x_1,x_2)$. So $x_1=\overline{g_1}$ and $x_2=\overline{g_2}$. Therefore $g_1=\alpha x_1$ and $g_2=\beta x_2+\gamma x_1$, where $\alpha,\beta,\gamma\in F$ and $\alpha,\beta\neq 0$, and then $\varphi$ is tame.\ \ \ $\square$

\end{document}